\documentclass[a4paper,11pt]{my_IJMS}

\usepackage{colortbl}
\usepackage{graphicx}
\newcommand{\menor} {\ {\raise-.5ex\hbox{$\buildrel<\over\sim$}}\ }
\newcommand{\maior} {\ {\raise-.5ex\hbox{$\buildrel>\over\sim$}}\ }

\theoremstyle{remark}
\newtheorem{example}{Example}[section]

\begin{document}


\title{A mu-differentiable Lagrange multiplier rule\footnote{Supported by {\it Centre for Research on Optimization and Control} (CEOC) from the ``Funda\c{c}\~{a}o para a Ci\^{e}ncia e a Tecnologia'' (FCT), cofinanced by the European Community Fund FEDER/POCI 2010. Accepted (22/June/2008) for
International Journal of Mathematics and Statistics (IJMS),
Vol.~4, No.~S09, Spring 2009 (in press).}}

\author{\textbf{Ricardo Almeida and Delfim F. M. Torres}}

\date{Department of Mathematics\\
University of Aveiro\\
3810-193 Aveiro, Portugal\\[0.3cm]
\{ricardo.almeida, delfim\}@ua.pt\\}

\maketitle


\begin{abstract}
\noindent \emph{We present some properties of the gradient of a mu-differentiable function. The Method of Lagrange Multipliers for mu-differentiable functions is then exemplified.}

\medskip


\medskip


\end{abstract}


\section{Introduction}

In \cite{AlmeidaTorres} we introduce a new kind of differentiation, what we call mu-differentiability, and we prove necessary and sufficient conditions for the existence of extrema points. For the necessary background on Nonstandard Analysis and for notation, we refer the reader to \cite{AlmeidaTorres} and references therein. Here we just recall the necessary results.

\begin{definition} \cite{AlmeidaTorres} Given an internal function $f:{^*\mathbb R}^n\to {^*\mathbb R}$, we say that $\alpha \in \mathbb R^n$ is a \emph{local m-minimum} of $f$ if
$$f(x) \maior f(\alpha) \, \mbox{ for all } \, x \in {^*B_r(\alpha)},$$
where $r\in \mathbb R$ is a positive real number.
\end{definition}

The crucial fact is that there exists a relationship between m-minimums and minimums:

\begin{lemma}\label{ponte2} \cite{AlmeidaTorres} If $f:{^*\mathbb R}^n\to {^*\mathbb R}$ is mu-differentiable, then
$$\alpha \mbox{ is a m-minimum of } f \mbox{ if and only if } \alpha \mbox{ is a minimum of } st(f).$$
\end{lemma}

With this lemma, and using the fact that

\begin{equation}\label{eq3}st \left( \left. \frac{\partial f}{\partial x_i} \right|_{\alpha} \right) = \left. \frac{\partial st(f)}{\partial x_i} \right|_{\alpha} \, \mbox{ for } \, i \in \{1,\ldots,n \},\end{equation}

it follows:

\begin{theorem} \cite{AlmeidaTorres} If $f:{^*\mathbb R}^n\to {^*\mathbb R}$ is a mu-differentiable function and $\alpha$ is a m-minimum of $f$, then
$$\left. \frac{\partial f}{\partial x_i} \right|_{\alpha}
\approx 0 , \mbox{ for every }\, i=1,\ldots,n.$$
\end{theorem}

In this paper we develop further the theory initiated in \cite{AlmeidaTorres}, proving some properties of the gradient vector (section \ref{gradient}) and a Method of Lagrange Multipliers (section \ref{MLM}). Illustrative examples show the analogy with the classical case.


\section{The Gradient Vector}\label{gradient}

In the sequel $f$ denotes an internal mu-differentiable
function from ${^*\mathbb R}^n$ to ${^*\mathbb R}$.

\begin{definition}
\label{def:grad}
A \emph{gradient vector} of $f$ at $x\in ns({^*\mathbb R}^n)$ is defined by
$$\nabla f(x) := \left( \left. \frac{\partial f}{\partial x_1} \right|_{x}, \ldots, \left. \frac{\partial f}{\partial x_n} \right|_{x} \right)$$
where
$$\left. \frac{\partial f}{\partial x_i} \right|_{x} \approx \frac{f(x_1,\ldots, x_{i-1}, x_i+\epsilon, x_{i+1}, \ldots,x_n)-f(x_1,\ldots,x_n)}{\epsilon}$$
and $\epsilon$ is an infinitesimal satisfying
$|\epsilon|>\delta_f$.
\end{definition}

\begin{remark}
The positive infinitesimal $\delta_f$
that appears in Definition~\ref{def:grad} is given
by the m-differentiability of $f$ (\textrm{cf.}
\cite{AlmeidaTorres}).
\end{remark}

\begin{remark}
Observe that
$$\left. \frac{\partial f}{\partial x_i}\right|_x \approx Df_x(e_i),$$
where $e_i = (0,\ldots, 0,1,0,\ldots,0)$ denotes the
$i$th canonical vector, and $Df_x$ denotes the derivative operator of $f$ at $x$.
\end{remark}

\begin{theorem} If $x,y \in ns({^*\mathbb R}^n)$ and $x\approx y$,
then
$$\left. \frac{\partial f}{\partial x_i}\right|_x \approx
\left. \frac{\partial f}{\partial x_i}\right|_y, \, \,
i=1,\ldots,n \, ,$$ \textit{i.e.}, $\nabla f(x) \approx \nabla f(y)$.
\end{theorem}

\begin{proof}
Simply observe that $Df_x(e_i)\approx Df_y(e_i)$.
\end{proof}

\begin{theorem} If $u\in{^*\mathbb R}^n$ is a finite vector, then
$$\forall x \in ns({^*\mathbb R}^n) \hspace{.5cm} Df_x(u) \approx \nabla f(x) \cdot u.$$
\end{theorem}

\begin{proof}
Since $st(f)$ is a $C^1$ function, if follows that for
any $v \in \mathbb R^n$
$$D st(f)_{st(x)}(v)=\nabla st(f)(st(x))\cdot v.$$
By the Transfer Principle of Nonstandard Analysis, it still holds for $u\in{^*\mathbb
R}^n$. On the other hand,
\begin{enumerate}
\item $D st(f)_{st(x)}(v)=st(Df_{st(x)})(u)\approx Df_x(u)$,
\item $\nabla st(f)(st(x))=st(\nabla f(st(x)))\approx \nabla
f(x)$,
\end{enumerate}
which proves the desired.
\end{proof}

We point out that, in opposite to classical functions, if  $\nabla f(x)$ is a gradient vector of $f$ at $x$, then $\nabla f(x)+\Omega$, where $\Omega\in{^*\mathbb R}^n$ is an infinitesimal vector, is also a gradient vector at $x$. Conversely, if $\nabla f(x)$ and $\nabla^1 f(x)$ are two gradient vectors, then $\nabla f(x)-\nabla^1 f(x)\approx 0$.

From now on, when there is no danger of confusion, we simply write $\nabla f$ instead of $\nabla f(x)$.

\begin{example}
Let $f(x,y,z)=(1+\epsilon)xy^2-\delta z$, with $(x,y,z) \in {^*\mathbb R}^3$, and $\epsilon$ and $\delta$ be two infinitesimal numbers. Given an infinitesimal $\theta$,
$$\begin{array}{rcl}
\displaystyle \frac{(1+\epsilon)(x+\theta)y^2-\delta z-((1+\epsilon)xy^2-\delta z)}{\theta}&=&(1+\epsilon)y^2,\\
&&\\
\displaystyle \frac{(1+\epsilon)x(y+\theta)^2-\delta z-((1+\epsilon)xy^2-\delta z)}{\theta}&=&2(1+\epsilon)xy+\theta (1+\epsilon)x,\\
&&\\
\displaystyle \frac{(1+\epsilon)xy^2-\delta (z+\theta)-((1+\epsilon)xy^2-\delta z)}{\theta}&=&-\delta,\\
\end{array}$$
and we can choose
$$\displaystyle \frac{\partial f}{\partial x}=(1+\epsilon)y^2, \, \, \displaystyle \frac{\partial f}{\partial y}=2(1+\epsilon)xy, \, \, \displaystyle \frac{\partial f}{\partial z}=-\delta.$$
\end{example}

\begin{theorem} If $f$ and $g$ are mu-differentiable and $k \in fin({^*\mathbb R})$, then
$$\nabla(kf)=k\nabla f, \quad \nabla(f+g)=\nabla f + \nabla g, \quad \mbox{ and } \quad \nabla(fg)=f\nabla g + g\nabla f.$$
\end{theorem}

\begin{proof} We prove only the last equality. Fix an infinitesimal number $\epsilon$ such that $|\epsilon|>\delta_f$. Then,
$$\frac{\partial (fg)}{\partial x_i} \approx \frac{(fg)(x_1,\ldots, x_{i-1}, x_i+\epsilon, x_{i+1}, \ldots,x_n)-(fg)(x_1,\ldots,x_n)}{\epsilon}$$
$$=f(x_1,\ldots,x_n) \frac{g(x_1,\ldots, x_{i-1}, x_i+\epsilon, x_{i+1}, \ldots,x_n)-g(x_1,\ldots,x_n)}{\epsilon}$$
$$+g(x_1,\ldots, x_{i-1}, x_i+\epsilon, x_{i+1}, \ldots,x_n)\frac{f(x_1,\ldots, x_{i-1}, x_i+\epsilon, x_{i+1}, \ldots,x_n)-f(x_1,\ldots,x_n)}{\epsilon}$$
$$\approx f(x)\frac{\partial g}{\partial x_i}+g(x)\frac{\partial f}{\partial x_i}$$
by the continuity of $g$.
\end{proof}

\begin{definition}
\label{def:m-critical-point}
We say that $x$ is a \emph{m-critical point} of $f$ if $\nabla f(x)\approx 0$.
\end{definition}

The following lemma is an immediate consequence of (\ref{eq3}) and Definition~\ref{def:m-critical-point}.

\begin{lemma} A point $x$ is a m-critical point of $f$ if and only if $st(x)$ is a critical point of $st(f)$.
\end{lemma}


\section{The Method of Lagrange Multipliers}\label{MLM}

Let $f:{^*\mathbb R}^n \to {^*\mathbb R}$ and $g_j:{^*\mathbb R}^n
\to {^*\mathbb R}$, $j=1,\ldots,m$ ($m \in \mathbb N$, $m<n$), denote
internal mu-differentiable functions. We address the problem of finding
m-minimums or m-maximums of $f$, subject to the conditions
$g_j(x)\approx0$, for all $j$. The constraints $g_j(x)\approx0$,
$j=1,\ldots,m$, are called \emph{side conditions}. Lagrange solved
this problem (for standard differentiable functions), introducing
new variables, $\lambda_1,\ldots,\lambda_m$, and forming the augmented function
$$F(x,\lambda_1,\ldots,\lambda_m)=f(x)+\sum_{j=1}^m \lambda_j g_j(x), \quad x \in \mathbb R^n.$$
Roughly speaking, Lagrange proved that the problem of finding the critical points of $f$, satisfying the
conditions $g_j(x)=0$, is equivalent to find the critical points
of $F$. We present here a method
to determine critical points for internal functions with
side conditions, based on the \emph{Method of Lagrange
Multipliers}. Similarly to the classical setting, define
$$\begin{array}{lcll}
F: & {^*\mathbb R}^{n+m} & \to & {^*\mathbb R}\\
  & (x_1,\ldots,x_n,\lambda_1,\ldots,\lambda_m) &
  \mapsto &  f(x_1,\ldots,x_n)+\displaystyle \sum_{j=1}^m \lambda_j g_j(x_1,\ldots,x_n).\\
\end{array}$$

If we let $g:=(g_1,\ldots,g_m)$ and
$\lambda:=(\lambda_1,\ldots,\lambda_m)$, we can simply write
\begin{equation}
\label{eq:def:F}
F(x,\lambda)=f(x)+\lambda \cdot g(x) \, .
\end{equation}

\begin{theorem}\label{LNF}[Lagrange rule in normal form with one constraint] Let $f:{^*\mathbb R}^n\to{^*\mathbb R}$ and $g:{^*\mathbb R}^n\to{^*\mathbb R}$ be two mu-differentiable
functions, and $\alpha$ a m-minimum of $f$ such that $g(\alpha)\approx 0$ and
$\nabla g(\alpha) \not\approx 0$. Then, there exists a finite $\lambda
\in {^*\mathbb R}$ such that
$$\nabla f(\alpha)+\lambda \nabla g(\alpha)\approx 0.$$
\end{theorem}

\begin{proof} Since $st(f)$ and $st(g)$ are functions of class
$C^1$, $\alpha$ is a minimum of $st(f)$, $st(g)(\alpha)=0$ and $\nabla st(g)(\alpha) \not=0$. It follows
(see, \textit{e.g.}, \cite[p.~148]{cheney}) that
$$\exists \lambda \in \mathbb{R} \hspace{.5cm} \nabla st(f)(\alpha)+\lambda \nabla st(g)(\alpha)=0.$$
Hence,
$$\nabla f(\alpha)+\lambda \nabla g(\alpha)\approx0.$$
\end{proof}

\begin{remark}
Suppose that we are in the conditions of Theorem \ref{LNF}. Then,
there exists some $\lambda_1\in fin({^*\mathbb R})$ such that
$$\nabla f(\alpha)+\lambda_1 \, \nabla g(\alpha)\approx0,$$
\textit{i.e.},
$$\left. \frac{\partial f}{ \partial x_i} \right|_{\alpha}+
 \lambda_1 \left. \frac{\partial g}{ \partial x_i} \right|_{\alpha}\approx0, \quad \, i=1, \ldots, n.$$
Using the notation \eqref{eq:def:F}, if $\alpha$
is a m-minimum of $f$ and $g(\alpha)\approx 0$, then
\begin{equation}\label{eq1}\left\{
\begin{array}{ll}
\left. \frac{\partial F}{ \partial x_i}
\right|_{(\alpha,\lambda_1)}=\left. \frac{\partial f}{\partial
x_i}\right|_{\alpha}+
 \lambda_1 \left. \frac{\partial g}{ \partial x_i} \right|_{\alpha} \approx 0, & i=1, \ldots,
 n \, ,\\
\left. \frac{\partial F}{ \partial \lambda}
\right|_{(\alpha,\lambda_1)}\approx g(\alpha)\approx0 \, . &\\
\end{array}
\right.
\end{equation}
Consequently, the m-critical points are solutions of the system
$$
\frac{\partial F}{ \partial x_i}\approx 0, \, i=1, \ldots,
n, \mbox{ and } \frac{\partial F}{ \partial \lambda}\approx0 \, ,
$$
\textit{i.e.}, $\nabla F \approx 0$.
\end{remark}

\begin{example}
Let $f(x,y,z)=xyz+\epsilon$, with $\epsilon\approx0$,
and consider the constraint $g(x,y,z)=x^2+2(y+\delta)^2+3z^2-1$,
with $\delta\approx0$. In this case, we define
$$F(x,y,z,\lambda):=xyz+\epsilon+\lambda
(x^2+2(y+\delta)^2+3z^2-1).$$ The system (\ref{eq1}) takes the form
$$\left\{
\begin{array}{ll}
yz+2\lambda x\approx0\\
xz+4\lambda(y+\delta)\approx0\\
xy+6\lambda z\approx0\\
x^2+2(y+\delta)^2+3z^2-1\approx 0 \, .
\end{array}
\right.$$
Since
$$xyz \approx -2\lambda x^2 \approx -4\lambda y(y+\delta) \approx -6\lambda z^2,$$
if $\lambda\not\approx0$, the solution is
$$x^2\approx\frac13, \, \, y^2\approx\frac16 \mbox{ and } z^2\approx\frac19 \, ;$$
if $\lambda \approx 0$, then
$$\left( 0,0,\pm \frac{1}{\sqrt 3} \right), \, \left( 0,\pm \frac{1}{\sqrt 2},0 \right)
 \mbox{ and } \left( \pm1,0,0\right)$$
are solutions. Observe that
$$\nabla g= (2x,4(y+\delta),6z) \approx (0,0,0) \mbox{ if and only if } (x,y,z)\approx(0,0,0).$$
One easily checks that
$$f\left(  \frac{1}{\sqrt 3}, \frac{1}{\sqrt 6}, \frac{1}{3} \right)=\frac{1}{3\sqrt{18}}+\epsilon
 \mbox{ is the m-maximum and }$$
$$f\left( - \frac{1}{\sqrt 3}, \frac{1}{\sqrt 6}, \frac{1}{3} \right)=-\frac{1}{3\sqrt{18}}+\epsilon
\mbox{ is the m-minimum}$$
of $f$ subject to the constraint $g$.
\end{example}

We now prove a more general Lagrange rule, admitting possibility of abnormal critical points ($\mu = 0$) and multiple constraints.

\begin{theorem}[Lagrange rule]
\label{thm:LR}
Let $f,g_1,\ldots,g_m$ be mu-differentiable
functions on ${^*\mathbb R}^n$. Let $\alpha$ be a m-minimum of $f$ satisfying
$$g_1(\alpha)\approx \ldots \approx g_m(\alpha)\approx 0.$$
Then, there exist finite hyper-reals
$\mu,\lambda_1,\ldots,\lambda_m \in {^*\mathbb R}$, not all
infinitesimals, such that
\begin{equation*}
\mu \nabla f(\alpha)+\lambda_1 \nabla g_1(\alpha)+\ldots+\lambda_m \nabla g_m(\alpha)\approx 0.
\end{equation*}
\end{theorem}

\begin{remark}
Defining $F(x,\mu,\lambda):=\mu f(x)+\lambda \cdot g(x)$, the necessary optimality condition given by Theorem~\ref{thm:LR} can be written as $\partial F / \partial x \approx \partial F / \partial \lambda \approx 0$.
\end{remark}

\begin{proof} First observe that $st(f),st(g_1),\ldots,st(g_m)$ are all
functions of class $C^1$, $\nabla st(f)(\alpha)=st(\nabla f)(\alpha)$ and
$\nabla st(g_j)(\alpha)=st(\nabla g_j)(\alpha)$, for $j=1,\ldots,m$.
Furthermore, since $\alpha$ is a minimum of $st(f)$ and
$$st(g_1)(\alpha)=\ldots =st(g_m)(\alpha)=0,$$
there exist reals $\mu,\lambda_1,\ldots,\lambda_m$, not all zero,
such that
$$\mu \, \nabla st(f)(\alpha)+\lambda_1 \nabla st(g_1)(\alpha)+\ldots+\lambda_m \nabla st(g_m)(\alpha)=0$$
(see, \textit{e.g.}, \cite[p.~148]{cheney}). Consequently,
\begin{equation}\label{eq10}
\mu \, st( \nabla f)(\alpha)+\lambda_1 st( \nabla
g_1)(\alpha)+\ldots+\lambda_m st( \nabla g_m)(\alpha)=0.
\end{equation}
On the other hand, we have
$$\mu \, st( \nabla f)(\alpha)=\mu \, st(\nabla f(\alpha))\approx \mu \nabla f(\alpha).$$
Analogously, for each $j=1,\ldots,m$,
$$\lambda_j st( \nabla g_j)(\alpha)\approx \lambda_j \nabla g_j(\alpha).$$
Substituting on equation (\ref{eq10}) the previous relations, one
proves the desired result.
\end{proof}

\begin{example}
Let $f(x,y,z)=z^2/2-(x+\epsilon)y$, with $\epsilon\approx0$,
be the function to be extremized,
and $g_1(x,y,z)=x^2+y-1$ and $g_2(x,y,z)=x+z-1+\delta$, with $\delta\approx0$, be the constraints. Then, the augmented function is
$$
F(x,y,z,\mu,\lambda_1,\lambda_2)
=\mu \left[z^2/2-(x+\epsilon)y\right]
+\lambda_1(x^2+y-1)+\lambda_2(x+z-1+\delta).$$
To find the local extrema of $f$, subject to the conditions $g_1\approx0$ and $g_2\approx0$, we form the system
\begin{equation}
\label{eq:nc:ex:t2}
\left\{\begin{array}{l}
-\mu y+2\lambda_1x+\lambda_2\approx0\\
- \mu (x+\epsilon)+\lambda_1\approx0\\
\mu z+\lambda_2\approx0\\
x^2+y-1\approx0\\
x+z-1+\delta\approx0 \\
\end{array}\right.
\end{equation}
of necessary optimality conditions. Assume $\mu \approx 0$
(abnormal case). Then, the first two equations in
\eqref{eq:nc:ex:t2} imply immediately that $\lambda_1 \approx
\lambda_2 \approx 0$. This is not a possibility by
Theorem~\ref{thm:LR}. We conclude that $\mu \not\approx 0$. The
solutions of \eqref{eq:nc:ex:t2} are then infinitely close to the
vectors
$$(-1,0,2) \quad \mbox{and} \quad (2/3,5/9,1/3).$$
Hence, if $f$ has any m-extrema under the given constraints, then
they must occur at either $(-1,0,2)$ or $(2/3,5/9,1/3)$.
\end{example}



\end{document}